\documentclass[12pt]{article}
\makeatletter
\newfont{\footsc}{cmcsc10 at 8truept}
\newfont{\footbf}{cmbx10 at 8truept}
\newfont{\footrm}{cmr10 at 10truept}
\renewcommand{\ps@plain}{%
\renewcommand{\@oddfoot}{\footsc
  the electronic journal of combinatorics {\footbf 11(2)} (2004--2005),
  \#A4\hfil\footrm\thepage}}
\makeatother
\newlength{\squarewidth} 
\newlength{\diagwidth} 
\font\shmt=chess20
\def\sqskip{\hspace*\squarewidth}
\def\ltod#1{
  \if#1KJ\else
  \if#1kj\else
  \if#1QL\else
  \if#1ql\else
  \if#1RS\else
  \if#1rs\else
  \if#1BA\else
  \if#1ba\else
  \if#1NM\else
  \if#1nm\else
  \if#1SM\else
  \if#1sm\else
  \if#1PO\else
  \if#1po\else
  Z
  \fi\fi\fi\fi\fi\fi\fi\fi\fi\fi\fi\fi\fi\fi
  }
\def\llhack#1{                
  \if#10\sqskip\else          
  \if#1.\sqskip\else
  \if#1SN\else\if#1sn\else#1\fi\fi\fi\fi}
\def\evenrank#1#2#3#4#5#6#7#8{
  \llhack#1\ltod#2\llhack#3\ltod#4\llhack#5\ltod#6\llhack#7\ltod#8}
\def\oddrank#1#2#3#4#5#6#7#8{
  \ltod#1\llhack#2\ltod#3\llhack#4\ltod#5\llhack#6\ltod#7\llhack#8}
\def\Board#1{\def\Heading{{\bf#1}}\board}
\def\Headingheight{1.2}
\def\board#1#2#3#4#5#6#7#8#9{\noindent
  \setlength{\unitlength}{1\squarewidth}
  \begin{picture}(9.5,10.8)
  {\shmt
  \put(.75,9.6){\makebox(8,\Headingheight)[c]{\Heading}}
  \put(.75,8.52){\evenrank #1}
  \put(.75,7.52){\oddrank  #2}
  \put(.75,6.52){\evenrank #3}
  \put(.75,5.52){\oddrank  #4}
  \put(.75,4.52){\evenrank #5}
  \put(.75,3.52){\oddrank  #6}
  \put(.75,2.52){\evenrank #7}
  \put(.75,1.52){\oddrank  #8}}
  \thinlines
  \put(.75,1.5){\framebox(8,8){}}
  \thicklines
  \put(.85,1.35){\line(1,0){8.07}}\put(8.9,9.4){\line(0,-1){8.07}}
  \put(.75,0.5){\makebox(8,1)[b]{\small #9}}
  \end{picture}
  }
\usepackage{amssymb,bm}
\def\0{^{\phantom0}}
\def\ns{$\!\!$} 
\def\BV{Bonsdorff\mbox{\kern.06em}-\mbox{\kern-.16em}V\"ais\"anen}
\title{
New directions in enumerative chess problems \\
 \normalsize\em to Richard Stanley on the occasion of his 60th birthday
}

\author{Noam D. Elkies
 \\ \normalsize  Department of Mathematics \vspace*{-0.5ex}
 \\ \normalsize  Harvard University \vspace*{-0.5ex}
 \\ \normalsize  Cambridge, MA 02138
 \\ \normalsize  \texttt{elkies@math.harvard.edu}
}

\date{\small 
Submitted: June 30, 2005;  Accepted: August 1, 2005;
Published: August 24, 2005}

\begin{document}
\maketitle

\begin{abstract}
Normally a chess problem must have a unique solution,
and is deemed unsound even if there are alternatives
that differ only in the order in which the same moves are played.
In an enumerative chess problem,
the set of moves in the solution is (usually) unique
but the order is not, and the task is to count the feasible permutations
via an isomorphic problem in enumerative combinatorics.
Almost all enumerative chess problems have been ``series-movers''$\!$,
in which one side plays an uninterrupted series of moves,
unanswered except possibly for one move by the opponent at the end.
This can be convenient for setting up enumeration problems,
but we show that other problem genres also lend themselves
to composing enumerative problems.  Some of the resulting enumerations
cannot be shown (or have not yet been shown) in series-movers.

This article is based on a presentation given at the banquet
in honor of Richard Stanley's 60th birthday,
and is dedicated to Stanley on this occasion.
\end{abstract}

{\large\bf 1 Motivation and overview}

Normally a chess problem must have a unique solution,
and is deemed unsound even if there are alternatives
that differ only in the order in which the same moves are played.
In an {\em enumerative chess problem},
the set of moves in the solution is (usually) unique
but the order is not, and the task is to count the feasible permutations
via an isomorphic problem in enumerative combinatorics.
Quite a few such problems have been composed and published
since about 1980 (see for instance \cite{Puusa, Stanley:ST}).
As Stanley notes in \cite{Stanley:ST}, almost all such problems
have been of a special type known as ``series-movers''$\!$.
In this article we give examples
showing how several other kinds of problems,
including the familiar ``mate in $n$ moves''$\!$,
can be used in the construction of enumerative chess problems.
We also extend the range of enumeration problems shown.  For instance,
we give a problem whose number of solutions in $n$ moves
is the \hbox{$n$\/-th} Fibonacci number, and another problem
that has exactly $10^6$ solutions.

This article is organized as follows.
After the above introductory paragraph and the following
Acknowledgements, we give some general discussion of
enumerative chess problems and of how a problem might
meaningfully combine mathematical content and chess interest.
We then introduce some more specific considerations
with two actual problems: one of the earliest enumerative
chess problems, by Bonsdorff and V\"ais\"anen,
and a recently composed problem by Richard Stanley.
We then challenge the reader with ten further problems:
another one by Stanley, and nine that we composed
and are published here for the first time.
We conclude by explaining the solution and mathematical context
for each of those ten problems.

{\bf Acknowledgements.}
This article is based on a presentation titled
``How do I mate thee?  Let me count the ways''
that I gave the banquet of the conference
in honor of Richard Stanley's 60th birthday;
the article is dedicated to him on this occasion.
I thank Richard for introducing me to queue problems
and to many other kinds of mathematical chess problems.
Thanks too to Tim Chow, one of the organizers of the conference,
for soliciting the presentation and proofreading a draft
of this article; to Tim Chow and Bruce Sagan,
for encouraging me to write it up for the present Festschrift;
and to the referee, for carefully reading the manuscript
and in particular for finding a flaw in the first version
of Problem~8.

This paper was typeset in \LaTeX, using Piet Tutelaers' chess font
for the diagrams.  Several of the problems were checked
for soundness with Popeye, a program created by Elmar Bartel,
Norbert Geissler, and Torsten Linss to solve chess problems.
The research was made possible in part
by funding from the National Science Foundation.

{\bf General considerations.}
All enumerative chess problems of the kind we are considering
lead to questions of the form ``in how many ways can one get
{}from position~X to position~Y in $n$ moves?''$\!$.\footnote{
  It would be interesting to have enumerative chess problems
  not of this form, which would thus connect
  chess and enumerative combinatorics in an essentially new way.
  To be sure, there are other known types of enumerative problems
  using the chessboard or chess pieces, but
  these are all chess {\em puzzles} rather than chess {\em problems},
  in that they use the board or pieces without reference
  to the game of chess.  The most familiar examples
  are the enumeration of solutions to the Eight Queens problem
  (combinatorially, maximal Queen \hbox{co-cliques}
  on the $8\times 8$ board) and of Knight's tours,
  and their generalizations to other rectangular board sizes.
  Of even greater mathematical significance are ``Rook numbers''
  (which count Rook co-cliques of size~$n$\/ on a given subset
  of an $N\times N$\/ board, see \cite[p.~71ff.]{Stanley:EC1})
  and the enumeration of tilings by dominos (a.k.a.\ matchings)
  of the board and various subsets
  (as in \cite[pp.~273--4 and 291--2, Ex.36]{Stanley:EC1};
  see also \cite{KELP}).
  }
But in general they are not explicitly formulated in this way,
because this would be too trivial in several ways.
It would be too easy for the composer
to pose an enumerative problem in this form;
it would be too easy for the solver to translate the problem
back to pure combinatorics;
and the problem would have so little chess content
that one could more properly regard it
as an enumerative combinatorics problem in a transparent chess disguise
than as an enumerative chess problem.  Instead, the composer
usually specifies only position~X, and requires that Y be
checkmate or stalemate.  (These are the most common goals
in chess problems, though one occasionally sees chess problems
with other goals such as double check or pawn promotion.)
The composer must then ensure that Y is the only such position
reachable within the stated number of moves, and the solver must first
find the target position~Y using the solver's knowledge or intuition
of chess before unraveling the problem's combinatorial structure.
This also means that one diagram suffices to specify the problem.
Another way to attain these goals is to exhibit only position~Y
and declare that X is the initial position where all 16 men
of one or both sides stand at the beginning of a chess game.
Most of the new problems in this article are of this type,
known in the chess problem literature as ``proof games''
or ``help-games'' (we explain this terminology later).

{\large\bf Two illustrative problems: \BV\ and Stanley}

\vspace*{1ex}
\centerline{
\Board{A: \BV, 1983}
{k.......}
{........}
{pPK.....}
{p.......}
{........}
{........}
{........}
{........}
{Series helpmate in 14.  How many solutions?}
\qquad
\Board{B: Richard Stanley, 2003}
{........}
{.b......}
{.rp..p.k}
{......N.}
{..Kp....}
{..b..p..}
{........}
{..R....N}
{Series helpmate in 7.  How many solutions?}
}

\vspace*{-1ex}

An early example of an enumerative chess problem is Diagram~A,
composed by Bonsdorff and V\"ais\"anen and published in 1983
in the Finnish problem periodical {\em Suomen Teht\"av\"aniekat}.
This problem, and Stanley's Diagram~B,
are examples of the ``series helpmate''$\!$,
an unorthodox genre of chess problems that is particularly well suited
to the construction of enumerative problems.
Black makes an uninterrupted series of moves,
at the end of which White has a (unique) mate in one.
The moves must be legal, and Black may not give check,
except possibly on the final move of the series
(in which case White's mating move must also parry the check).
Problems that require one side to make a series of moves
are known as ``series-movers''$\!$.  Series stipulations
appear regularly in the problem literature, though they are regarded
as unorthodox compared to stipulations in which White and Black
alternate moves as in normal chess-play.  Such alternation
is not a common element in enumerative combinatorics,
and most enumerative chess problems avoid it,
either explicitly by using a series stipulation,
or implicitly by ensuring that the combinatorial structure
involves only one player's moves.
This is the case for almost all problems in this article.
A notable exception is Diagram~3, where (as in \cite{CEF})
the combinatorial problem is chosen to be expressible
in terms of move alternation.  In one of the other problems,
both White and Black moves figure in the enumeration
but do not interact, so that the problem reduces to a pair
of series-movers.  Likewise, enumerative problems usually
do not involve struggle between Black and White:
indispensable though it is to the game of chess,
this struggle does not easily fit into a combinatorial problem.
Usually the stipulation simply requires both sides to cooperate,
or one side not to play at all, thus pre-empting any struggle.
Our Diagram~4 is presented as a ``mate in $n$'' problem,
which usually presupposes that Black strives to prevent this mate;
but here Black has no choice, so again there is no real struggle.
In Diagram~5, also a ``mate in~$n$'', Black again can do nothing
to hinder White, but does have some choices,
which the solver must account for.

\vspace*{1ex}
\centerline{
\Board{A {\rm (again)}: \BV, 1983}
{k.......}
{........}
{pPK.....}
{p.......}
{........}
{........}
{........}
{........}
{Series helpmate in 14.  How many solutions?}
\qquad
\Board{A$\!\bm'$: \BV, 1983}
{kb......}
{bP......}
{..K.....}
{........}
{........}
{........}
{........}
{........}
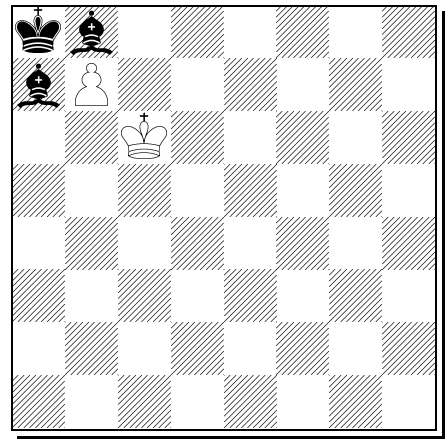
{The target position for Diagram~A}
}

\vspace*{-1ex}

In the \BV\ problem, Black has $14$ moves
to reach a position where White can give checkmate.
The only such checkmate reachable in as few as $14$ Black moves
is A$'$, after both pawns have promoted to Bishops
and moved to b8 and a7 via e5, blocking the King's escape
so that White's move b7 gives checkmate.
Thus the pawns/Bishops must travel along the following route:

\vspace*{1ex}
\centerline{
{\bf a6}---{\bf\underline{a5}}\/---a4---a3---a2---a1(B)---e5---{\sf b8}---{\sf\underline{a7}}
}

starting at a6 and a5, moving one space at a time, and ending at b8 and a7,
with the pawn that starts on a5 (and the Bishop it promotes to)
always in the lead.  Enumerative chess problems such as this, where
all the relevant chessmen move in one direction along a single path,
are known as ``queue problems'': the chessmen are imagined to be waiting
in a queue and must maintain their relative order.
The number of feasible move-orders is given by a known but nontrivial
formula, making such queues appropriate for an enumerative chess problem.
Here the queue contains just two units, which begin
at the first two squares of the path and end in its last two squares.
In this case the formula yields $C_n = (2n)!/(n!(n+1)!)$,
where $n$ is the number of moves played by each unit in the queue.
Hence the number of solutions of Diagram~A is
$\bm{C_7 = 14!/7!8! = 429}$.

The $C_n$ are the celebrated {\em Catalan numbers},
which enumerate a remarkable variety of combinatorial structures;
see \cite[pages 221--231]{Stanley:EC2}\footnote{
  Also available and updated online from Richard Stanley's website,
  see \cite{Stanley:CAT}.
  }
and Sequence~A000108 in \cite{Sloane:OEIS}.
In the setting of enumerative chess problems,
a particularly useful way to see that Diagram~1 has $C_7$ solutions
is to organize Black's moves as follows:

\vspace*{1ex}

\centerline{
\begin{tabular}{ccccccc}
\ns a4 \ns & \ns a3 \ns & \ns a2 \ns & \ns a1B \ns &
  \ns Be5 \ns & \ns Bb8 \ns & \ns Ba7 \ns \\
\ns a5 \ns & \ns a4 \ns & \ns a3 \ns & \ns a2 \ns &
  \ns a1B \ns & \ns Be5 \ns & \ns Bb8
\end{tabular}
}

The top (resp.\ bottom) row contains the lead (rear) pawn's moves;
a move order is feasible if and only if
each move occurs before any other move(s) appearing
in the quarter-plane extending down and to the right from it.
These constraints amount to a structure of a {\em poset}\/
({\em p}artially {\em o}rdered {\em set}) on the set of Black's moves,
with $x \prec y$ if and only if $x \neq y$
and move~$y$ appears in or below the row of~$x$,
and in or to the right of the column of~$x$.
In the problem, $x \prec y$ means $x$ must be played before~$y$,
and a solution amounts to a {\em linear extension}\/ of~$\prec$,
that is, a total order consistent with~$\prec$.
This kind of analysis applies to many enumerative chess problems.
There is no general formula for counting
linear extensions of an arbitrary partial order,
but in many important cases a nontrivial closed form is known.
For a queue problem such as Diagram~A, with two chessmen
that start next to each other at one end of the route
and finish next to each other at the other end, the poset is
the Young diagram corresponding to the partition $(n,n)$ of~$2n$,
and a linear extension corresponds to a standard Young tableau
of shape $(n,n)$.  Therefore the number $C_n$ of extensions
can be obtained from the hook-length formula.

The hook-length formula also answers any queue problem with
$k$\/ chessmen that start at the first $k$ squares
of the queue line, or equivalently end on the last $k$ squares.
Many such problems have been composed (see for instance \cite{Puusa}).
Even the special case of $k=2$ queues that lead to Catalan numbers
has appeared in several published problems besides the \BV\ problem
analyzed here.  One example is a V\"ais\"anen problem that appears
as Exercise~6.23 in \cite[p.232]{Stanley:EC2}.
Another is the problem cited as Diagram~0 in the next section.

\vspace*{-1ex}

\def\Headingheight{.75}
\centerline{
\Board{B {\rm (again)}: Richard Stanley, 2003}
{........}
{.b......}
{.rp..p.k}
{......N.}
{..Kp....}
{..b..p..}
{........}
{..R....N}
{Series helpmate in 7.  How many solutions?}
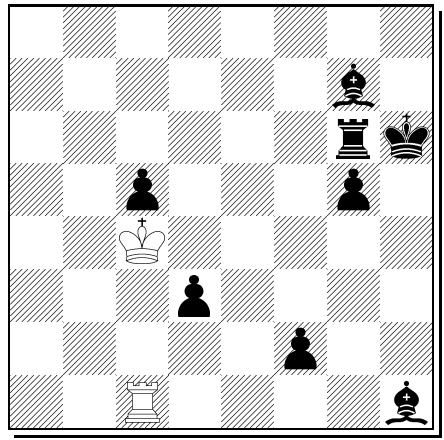
\quad\qquad
\Board{B$\bm'$: Richard Stanley, 2003}
{........}
{......b.}
{......rk}
{..p...p.}
{..K.....}
{...p....}
{.....p..}
{..R....b}
{The position after Black's series in Diagram~B}
}

\def\Headingheight{1}

\vspace*{-2ex}

\pagebreak

Diagram B is a problem by Stanley \cite[pp.7--8]{Stanley:ST}
that also leads to an enumeration of linear extensions
of a partial order, but one of a rather different flavor.
Black must play the four pawn moves c6, d3, f2, fxg5,
opening lines for Black's Rook and two Bishops to play
Rg6, Bg7, Bxh1 to reach position~B$'$, after which White mates with Rxh1.
In a feasible permutation, each Rook or Bishop move must be played
after its two line-opening pawn moves.  We write these constraints as

\vspace*{1ex}

\centerline{
f2 $<$ Bxh1 $>$ c5 $<$ Rg6 $>$ fxg5 $<$ Bg7 $>$ d3.
}

This means that in any feasible order of Black's moves, such as

\vspace*{1ex}

\centerline{
1~c5 \ 2~fxg5 \ 3~f2 \ 4~Rg6 \ 5~d3 \ 6~Bxh1 \ 7~Bg7,
}

the moves f2, Bxh1, \ldots, d3 must be numbered by integers
that constitute a permutation of $\{1,2,\ldots,7\}$
satisfying those inequalities (such as
$$
3 < 6 > 1 < 4 > 2 < 7 > 5
$$
in our example).  Therefore the solutions of Diagram~B
correspond bijectively with up-down permutations\footnote{
  Also known as ``zigzag permutations'' or,
  confusingly (because there is no parity condition),
  ``alternating permutations''$\!$.  In~\cite{Stanley:ST}
  Stanley uses the convention that zigzag permutations
  $(\sigma_1, \sigma_2, \sigma_3, \ldots)$
  must satisfy $\sigma_1 > \sigma_2 < \sigma_3 > \; < \cdots$
  rather than $\sigma_1 < \sigma_2 > \sigma_3 < \; > \cdots$,
  and thus replaces each $\sigma_i$ by $8-\sigma_i$
  before constructing the bijection.
  }
of order~$7$.
It is known that the number of up-down permutations of order~$n$
is the \hbox{$n$\/-th} {\em Euler number} $E_n$,
which may be defined by the generating function
$$
\sec x + \tan x = \sum_{n=0}^\infty E_n \frac{x^n}{n!}
$$
(see for instance \cite[Problem~5.7]{Stanley:EC2}, \cite{NDE:zeta},
and Sequence~A000111 in \cite{Sloane:OEIS}).  Therefore
Diagram~B has $\bm{E_7 = 272}$ solutions.

{\large\bf Some new enumerative chess problems}

Diagram 0, reproduced from \cite{Stanley:ST},
is a queue problem composed by Stanley for his guest lecture
at the author's seminar on Chess and Mathematics
for Harvard freshmen.  This problem more than doubles the length
of the pawn/Bishop queue in the \BV\ problem
(Diagram~A of the Introduction), and is the longest such problem known.

The remaining problems in this article appear here for the first time.
Diagram~1 is to be reached cooperatively by White and Black from
the starting position in the minimal number of moves.  The moves,
however bizarre strategically, must obey all the rules of chess,
including those involving check: the ``cooperation'' does not extend
to letting the opponent put or leave the King in check,
nor to overlooking other illegal moves.
(This kind of problem is called a ``proof game''\footnote{
  Conventionally every chess problem, of whatever genre,
  must be reachable by a legal game from the initial position;
  a ``proof game'' ending in a given position thus proves
  that the position is legal.  In a proof game problem,
  the (usually minimal) length of the game is also specified,
  usually with the intention that this forces a unique
  and remarkable solution.  See for instance~\cite{SPG}
  for some good examples of what can be done in this genre.
  }
or, less confusing to a mathematician, a ``help-game''$\!$.)
The resulting enumeration problem has already been explained.

\centerline{
\Board{0: Richard Stanley, 2003}
{k.......}
{........}
{pPPp.P..}
{p.Pp....}
{...K...p}
{..P.P..p}
{.......p}
{......nr}
{Series helpmate in 34.  How many solutions?}
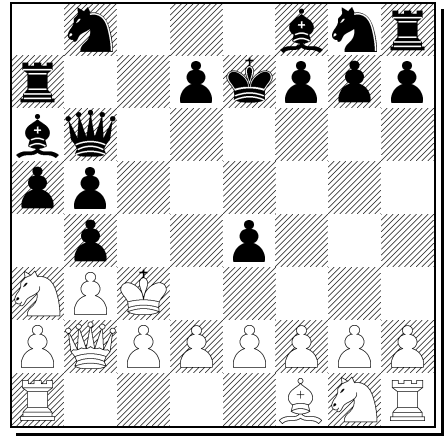
\qquad
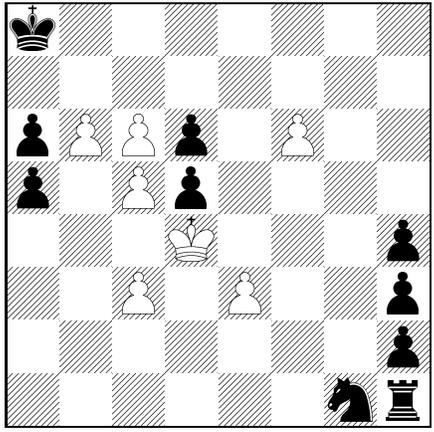
\Board{1: NDE, 2004}
{.n...bnr}
{r..pkppp}
{bq......}
{pp......}
{.p..p...}
{NPK.....}
{PQPPPPPP}
{R....BNR}
{How many shortest games?}
}

\vspace*{-1ex}

Diagram 2 leads to an enumeration problem that may be regarded
as a generalization of both of the types seen so far
(which gave Catalan and Euler numbers).
The stipulation is analogous to that of Diagram 1,
but involves only the White chessmen, which are
to reach the diagram from their initial array
in the least number of legal moves.
This is thus a kind of series-mover; when such problems are composed
to have a unique solution they are usually called ``series proof games''
or ``one-sided proof games'': with only one side playing, we cannot
speak of cooperation, and thus avoid the term ``help-game''$\!$.

\centerline{
\Board{2: NDE, 2004}
{........}
{........}
{........}
{.....PP.}
{...PPBQ.}
{...B....}
{PPP.N..P}
{RN..K..R}
{How many shortest sequences?}
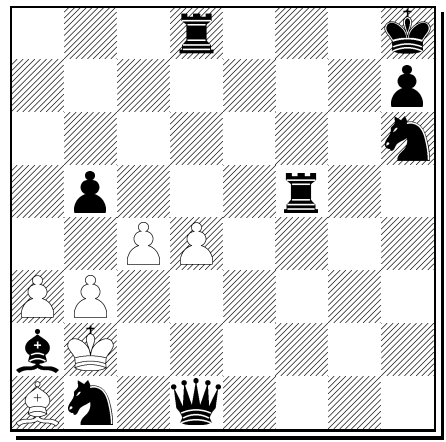
\qquad
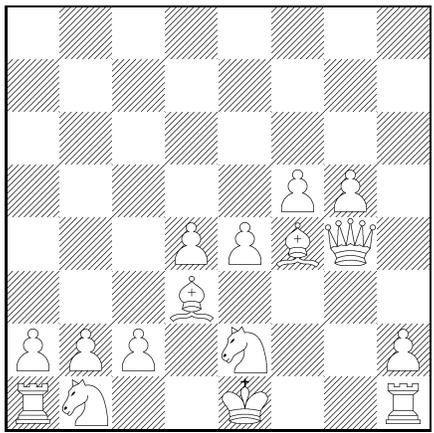
\Board{3: NDE, 2004}
{...r...k}
{.......p}
{.......n}
{.p...r..}
{..PP....}
{PP......}
{bK......}
{Bn.q....}
{Helpmate in 3.5\,.  How many solutions?}
}

\vspace*{-1ex}

Diagram 3 is a {\em helpmate}\/: Black and White cooperate
to get Black mated in the stipulated number of moves.
Here the move count of~3.5 means that White moves first,
and then Black helps White give checkmate on White's fourth move.
As with Diagram~1, all moves must be legal.
This leads to an enumerative problem recently introduced in~\cite{CEF},
and illustrated there by a help-stalemate.
Diagram~3 answers a challenge by Tim Chow,
one of the authors of~\cite{CEF}, to show this enumeration
in helpmate form.

Each of the next two problems shows an infinite sequence:
the \hbox{$n$-th} term of the sequence is the number of ways
White can force checkmate in exactly $n$ moves.
Both are much simpler than the number of pawns and pieces
might suggest, since most of these units are immobile
and serve only to restrict White's and Black's choices.
The intended answer to the second problem may be controversial.
In both problems, the solver should ignore the fifty-move rule
and the rule of triple repetition: in actual chess play
such rules are needed to oblige recalcitrant players to accept a draw
when neither side can force a win, but in most problems these rules
do not apply.

\centerline{
\Board{4: NDE, 2003--4}
{k.....b.}
{P....p..}
{KP...Pp.}
{PP...pP.}
{.....P..}
{.....p..}
{.....P..}
{......B.}
{Mate in (exactly) $n$: how many solutions?}
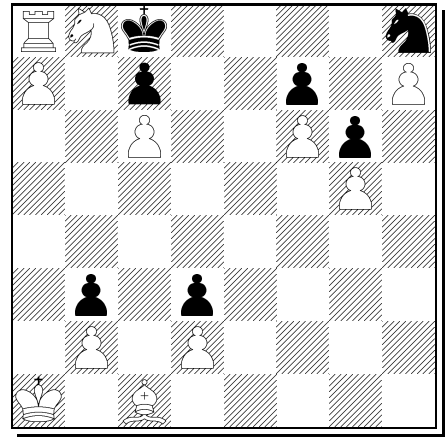
\qquad
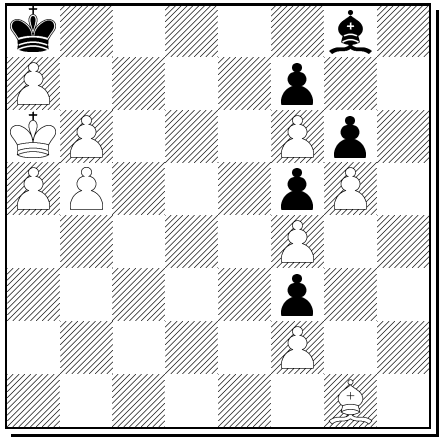
\Board{5: NDE, 2005}
{RNk....n}
{P.p..p.P}
{..P..Pp.}
{......P.}
{........}
{.p.p....}
{.P.P....}
{K.B.....}
{Mate in (exactly) $n$: how many solutions?}
}

\vspace*{-1ex}

Our final four problems were composed to attain a specific number
of numerological rather than mathematical interest.
The first two are series proof games.
Both were composed as New Year's greetings, and breach the convention
that all solutions must consist of the same set of moves.
Diagram~6 was also used as an ``entrance exam''
for the Chess and Mathematics seminar mentioned earlier,
see~\cite{NDE:2004}.
In the remaining two problems
we return to help-games with a unique move-set.
Diagram~8 was suggested by a helpmate by K.~Fabel
({\em Heidelberger Tageblatt}, 8.x.1960)
that has exactly $1000$ solutions in $5$ moves.
Diagram~9 was composed for Richard Stanley in honor of his 60th birthday
and was first presented at the banquet dinner of his birthday conference.

\centerline{
\Board{6: NDE, 1/2004}
{........}
{........}
{........}
{.K..N...}
{..P.....}
{..N..PP.}
{PP.PP..P}
{.RBQ.B.R}
{How many shortest sequences?}
\qquad
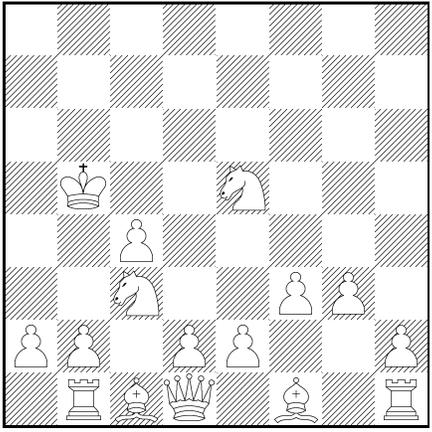
\Board{7: NDE, 12/2004}
{........}
{........}
{........}
{........}
{N...QP..}
{..PP...K}
{PPR.P.PP}
{..B..BNR}
{How many shortest sequences?}
}

\vspace*{1ex}

\centerline{
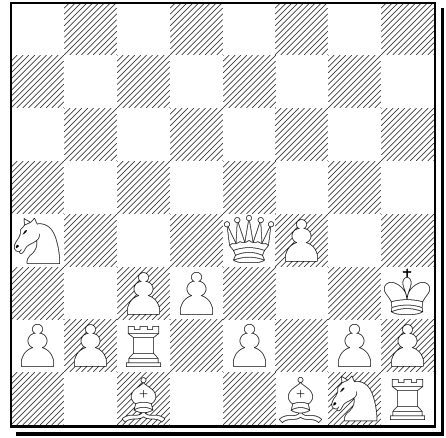
\Board{8: NDE, 2003}
{.n.k...r}
{..pp.ppp}
{...bbn..}
{.p..prq.}
{pP..P...}
{....B.PK}
{P.PPBP.P}
{RN.QRN..}
{How many shortest games?}
\qquad
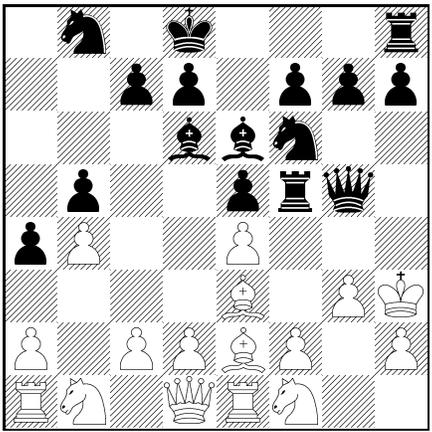
\Board{9: NDE, 2004, for RS-60}
{rn...bnr}
{ppp..ppp}
{...ppk..}
{...N..q.}
{......bN}
{.P......}
{PBPPPPPP}
{R..QKB.R}
{How many shortest games?}
}

\vspace*{2ex}

{\large\bf Solutions and Comments}

{\bf 0.} There are $\bm{C_{17} = 34!/17!18! = 129644790}$ solutions.
As in the Bonsdorff-V\"ais\"anen problem, two Black \hbox{a-pawns} promote
to Bishops and then travel to b8 and~a7 so that White's b7 is checkmate.
The pawns/Bishops travel along a unique path and never occupy
the same square at the same time.  But here the path is longer:
a6-a5-a4-a3-a2-a1B-b2-c1-d2-e1-g3-f4-h6-f8-e7-d8-c7-b8-a7.
Each pawn/Bishop makes $17$ moves, so the number of feasible
permutations of the $2\times 17 = 34$ moves is the 17th Catalan number.
The prohibition against checking before the final move of Black's
sequence is used extensively: Black's cluster around the h1 corner,
which serves only to block an alternative path through h4 and g5
(instead of g3-f4), is immobile because moving the Knight from g1
to either f3 or e2 would check the Kd4; likewise neither Black pawn
may promote to a Knight (which could reach a7 or b8 more quickly
than a Bishop), because the first move of a Knight from~a1 would check
White's King from b3 or~c2; and the Black Bishops must detour around
the White pawns at c3,e3,f6 because capturing any of those pawns
would again check the White King.
The c5 pawn blocks the line a7--d4 so that the \hbox{b6-pawn}
is not pinned by a Black Ba7 and may move to b7 to give checkmate.

{\bf 1.} There are $\bm{E_9 = 7936}$ solutions
of the minimal length of $10$ moves
(as usual a ``move'' comprises both a White and a Black turn).
Since White is in check from the Pb4, that pawn must have made
the last move, necessarily a capture from c5, and the only
missing White unit is the dark-square Bishop.
We quickly deduce that White must have played at least $10$ moves,
and could play exactly~$10$ only if they were
b3,Ba3,Bb4,Na3,Qb1,Kd1,Kc1,Kb2,Kc3,Qb2 in this order.
Black also needs $10$ moves to reach the diagram,
and there is only one set of $10$ moves that attains this:
a5,b5,c5,c5xb4,e5,e4,Ra7,Ba6,Qb6,Ke7.  We saw that c5xb4
must be played last, but there are many choices for the order
of the remaining $9$~moves.  By writing the constraints as

\vspace*{1ex}

\centerline{
Ra7 $>$ a5 $<$ Ba6 $>$ b5 $<$ Qb6 $>$ c5 $<$ Ke7 $>$ e5 $<$ e4,
}

we obtain a bijection between the feasible orders and the up-down
permutations of order~$9$, and find that there are $E_9 = 7936$
feasible orders, as claimed.  Therefore this problem answers
the challenge posed in \cite[p.8]{Stanley:ST}:
``Can the theme of [Diagram~B] be extended to
$E_8=1385$ or $E_9=7936$ solutions?''

While only Black moves figure directly in the enumeration,
White's $10$ moves are not entirely idle.  White not only
provides fodder and target for the final checkmate but also,
with the early Ba3, enforces the condition c5\,$<$\,Ke7:
playing Ke7 first would move the Black King into check.

Enumerative proof games leading to Catalan numbers and other queue
variants have also been constructed (by Andrew Buchanan and us),
but so far have not beaten the series-mover records.

%

{\bf 2.} There are $\bm{3850}$ solutions
of the minimal length of $10$~moves.
The number $3850$ arises as one-half the number of
standard Young tableaux (SYT's) associated with
the self-conjugate partition $\lambda=(5,3,2,1,1)$ of~$12$.

The unique set of $10$ moves that reach this position
{}from White's opening array, and the constraints on the order
in which they may be played, may be read off the following diagram:

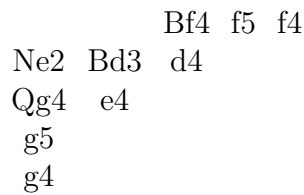
\vspace*{1ex}

\centerline{
\begin{tabular}{ccccc}
            &             & \ns Bf4 \ns & \ns f5 \ns & \ns f4 \ns \\
\ns Ne2 \ns & \ns Bd3 \ns & \ns d4 \ns \\
\ns Qg4 \ns & \ns e4 \ns \\
\ns g5 \ns \\
\ns g4 \ns
\end{tabular}
}

Each move must be played {\em after} any move or moves
to its right or below it.  Thus the number of feasible permutations
is the number of linear extensions of the partial order
given by the ``skew Young diagram'' $\lambda/\mu$,
where $\lambda=(5,3,2,1,1)$ as above and $\mu=(2)$.
This number, call it $P_{\lambda/\mu}$ (for any partitions
$\lambda,\mu$ such that $\mu_i \leq \lambda_i$ for all~$i$),
occurs in algebraic combinatorics,
notably in Pieri's rule (see for instance \cite[p.59]{FH}),
which asserts that $P_{\lambda/\mu}$ is the multiplicity
of the representation $V_\mu$ in the restriction of~$V_\lambda$
{}from~$S_{|\lambda|}$ to~$S_{|\mu|}$.
Each of our problems so far comes down to the evaluation
of $P_{\lambda/\mu}$ for some choice of $\lambda,\mu$;
for instance, Diagram~0 leads to $\lambda=(17,17)$ and $\mu=0$,
while Diagram~1 amounts to $\lambda=(5,5,4,3,2)$ and $\mu=(4,3,2,1)$.
In general there is no simple formula for $P_{\lambda/\mu}$,
but many special cases have nice enumerations.
For instance, if $|\mu|\leq1$ then
$P_{\lambda/\mu}$ is the number of SYT's of shape~$\lambda$
(since any such SYT must have~$1$ in the top left box),
and so can be computed by the hook-length formula.
In our case $|\mu|=2$.
The general formula for $|\mu|=2$ is more complicated,
but our $\lambda$ is self-conjugate,
so by symmetry exactly half the SYT's of shape~$\lambda$
have $1$ and~$2$ in the first two boxes of the top row.
These are exactly the tableaux that yield
feasible permutations of White's ten moves via
a standard skew Young tableau of shape $\lambda/\mu$.
Hence the number of feasible permutations is
$\frac12 \dim(V_\lambda) = 3850$, as claimed.\footnote{
  In the setting of Pieri's rule, $|\mu|=1$ makes $V_\mu$
  the trivial representation of the trivial group, so clearly
  $P_{\lambda/\mu}=\dim(V_\lambda)$.  For $|\mu|=2$,
  the general formula for $P_{\lambda/\mu}$ requires
  the computation of both $\dim(V_\lambda)$
  and the character of a simple transposition acting on $V_\lambda$;
  but when $\lambda$ is self-conjugate we have
  $V_\lambda \cong V_\lambda \otimes \varepsilon$,
  so the action of any odd permutation on~$V_\lambda$
  automatically has character~$0$.
  In this case $\frac12\dim(V_\lambda)$
  is the dimension of an irreducible representation $V'_\lambda$
  of the alternating group $A_{|\lambda|}$
  such that $V'_\lambda \oplus V'_\lambda$ is the restriction
  of $V_\lambda$ from $S_{|\lambda|}$ to $A_{|\lambda|}$.
  }

{\bf 3.} There are $\bm{2}$ solutions, namely\footnote{
  As in \cite[pp.14--15]{CEF}, we give the solutions
  with White beginning each move.  This reverses the usual
  problemists' convention for help-problems, which would write
  the first solution 1\ldots cxb5 2~Rxb5~b4 etc.,
  but is more natural for the great majority of chess-players
  who rarely encounter helpmates.
  }

\vspace*{1ex}

\centerline{\bf 1~cxb5 Rxb5 2~b4 Bg8 3~d5 Qa4 4~Kxb1\#,}
\centerline{\bf 1~b4 Qa4 2~cxb5 Bg8 3~Kxb1 Rxb5 4~d5\#.}

These are permutations of the same four White and three Black moves.
The ordering constraints on these moves are given by the following
diagram, in which Black's moves are shown in {\bf boldface}
as in~\cite[pp.14--15]{CEF}:

\vspace*{1ex}

\centerline{
\begin{tabular}{ccc}
             & \ns d5 \ns & \ns \ns {\bf Rxb5} \ns \ns \\
\ns Kxb1 \ns & \ns {\bf Bg8}  \ns & \ns  cxb5  \ns \\
\ns {\bf Qa4} \ns & \ns b4 \ns
\end{tabular}
}

Again each move must be played
after any move or moves to its right or below it.
Thus as in \cite{CEF} we seek extensions of this partial order
that also respect the alternation between White and Black moves.
This partial order corresponds to the $3\times 3$ Young diagram,
with the top left and bottom right corners removed ---
which does not change the number of linear extensions
because these are the minimum and maximum of the poset
and have the correct parity in the checkerboard coloring
of the $3\times 3$ diagram.
Hence the solutions correspond bijectively with
what are called in~\cite{CEF} ``chess tableaux'' of shape $(3,3,3)$.
Our solutions correspond to the two chess tableaux of this shape.

\pagebreak

{\bf 4.}
The number of mates in~$n$ is the $n$\/-th Fibonacci number $\bm{F_n}$,
that is $1,1,2,3,5,8,13,\ldots$ mates in $1,2,3,4,5,6,\ldots$ moves.
Until White mates by playing b7, only the Bishops can move:
Black's will shuttle between g8 and h7, and White's along the path
g1-h2-g3-h4 of length~$4$.  Thus the number of solutions equals
the number of walks of length $(n-1)$ on that path beginning
at the endpoint~g1, a number which is well-known (and easily shown)
to equal~$F_n$.
It is important that Black never has any choice,
as the next problem illustrates.

{\bf 5.}  We claim that there are $\bm{2^{F_n}}$ solutions.
We argue as follows.  Until White mates by playing the Knight
to a6 or d7, only the Kings can move:
White's will shuttle between a1 and b1,
and Black's along the path c8-d8-e8-f8 of length~$4$.
If White will checkmate on move~$n$ then Black has $F_n$ choices
for Black's $n-1$ moves, as we saw in the previous problem.
To completely specify how White will checkmate on the
\hbox{$n$\/-th} move, then, White must declare,
for each of Black's $F_n$ possible sequences,
a choice between Na6 and Nd7, and this can be done in $2^{F_n}$ ways!

It might appear that the White Rook and Knight are superfluous:
without them, the Kings are still confined to the same paths
until White promotes on a8 to a Queen or Rook, again with
two checkmating options against each of $F_n$ Black sequences.
But White could also promote to a Bishop,
or even to a Knight after Black plays Kf8.
White could then soon capture Black's pawns on d3 and b3,
freeing the immured King and Bishop, and eventually win,
producing innumerable extra solutions once $n$ is large enough.

{\em Exercise}\/:  Construct positions in which White has
$2^n$, $2^{2^{n-1}}$, or $2^{n2^{n-1}}$ ways to force checkmate
in $n$ moves.

{\bf 6.} There are $\bm{2004}$ sequences of the minimal length~$12$.
Each consists of the single move g3, the \hbox{$3$-move} sequence
c4,Nc3,Rb1, and one of the three \hbox{$8$-move} sequences
Nf3,Ne5,f3,Kf2,Ke3,Kd3(d4),Kc4(c5),Kb5.  The move~g3 may be played
at any point, and so contributes a factor of~$12$.
If the King goes through~c5 then the 3- and \hbox{$8$-move} sequences
are independent, and can be played in ${11 \choose 3}$ orders.
If the King goes through~c4 then the entire \hbox{$8$-move} sequence
must be played before the \hbox{$3$-move} sequence begins,
so there are only two possibilities,
depending on the choice of Kd3 or~Kd4.  Hence the total count is
$12 \bigl({11 \choose 3} + 2\bigr) = 2004$ as claimed.

{\bf 7.} There are $\bm{2005}$ sequences of the minimal length~$14$.
This and the next problem use the happy coincidence\footnote{
  Perhaps Scheherazade had only $14$ basic stories,
  and combined them in sets of~$4$ in all possible ways.
  }
${14 \choose 4} = 1001$.
Here White plays the \hbox{$4$-move} sequence f4,Kf2,Kg3,Kh3
and one of the five sequences
Nc3,Na4,c3,Qc2,Qe4,d3,Bd2(e3,f4,g5,h6),Rc1,Rc2,Bc1 of length~$10$.
If the Bishop goes to d2 or~e3, the sequences are independent,
and can be played in ${14 \choose 4}$ orders.
Otherwise the Bishop must return to c1 before White plays~f4,
so the entire \hbox{$10$-move} sequence
must be played before the \hbox{$4$-move} sequence begins.
Hence the total count is $2 {14 \choose 4} + 3 = 2005$.

{\bf 8.} Exactly $\bm{10^6}$.  White and Black play independently,
and each can reach the position in $1000$ ways in the minimal
number of moves ($14$ for White, $13$ for Black).  Curiously
neither the White nor the Black enumeration uses the factorization
$1000 = 10 \times 10 \times 10$ or $1000 = 2^3 \times 5^3$,
though the factor $2\times 4$ does figure in the Black enumeration.
White's $1000$ is ${14 \choose 4} - 1$: the $4$- and \hbox{$10$-move}
sequences b4,Bb2,Bd4,Be3 and \hbox{e4,Ne2,Ng3,Be2,0-0,Re1,Nf1,g3,Kg2,Kh3}
are independent except for the condition that e4 must precede Be3.
Black's $1000$ is $2 \times 4 \times \bigl({9 \choose 4} - 1\bigr)$.
Black starts a5,a4,Ra5,Rf5.
Then the \hbox{$4$-move} sequence b5,Bb7,Bd5,Be6
is independent of the \hbox{$5$-move} set e5,Qg5,Nf6,Kd8,Bd6
as long as e5 precedes Be6.  Of the latter $5$~moves,
e5 must come first, but then Black has choices: Nf6 and Kd8
can be played in either order after Qg5, and Bd6 can be interpolated
in any of $4$~spots, whence the factors of~$2$ and~$4$.
Note that there is no danger of White's King being in check on~g2
{}from Black's Bishop on b7 or d5,
because White's move e4 always precedes Kg2;
nor of the King's being in check on~h3,
because Black's Rf5 always precedes Be6.
In White's sequence it might seem that \hbox{0-0}
could be replaced by Kf1, but this would lead to a vicious circle:
g3 must then precede Kg2, which must precede Re1,
which must precede Nf1, which must precede g3 --- contradiction.
Hence White castles to get the King out of the way.

To my knowledge this is the first chess problem composed
to have exactly a million solutions.

{\bf 9.}  Here the target number was $\bm{60}$,
since this problem was composed for Stanley's 60th birthday.
Since it is easy to make $60$ the answer to an enumerative problem,
there was considerable scope for chess content.
In accordance with the title of the banquet presentation,
Black is checkmated in the diagram (by double check; Black could deal
with each of Bb2 and Nd5 separately, but not with both at once).
Each side needs only six moves to reach the diagram: Black by
the unique sequence d6,Bg4,e6,Qg5,Ke7,Kf6, White by
b3,Bb2,Nf3,Nh4,Nc3,Nd5 in some order.  But in fact White must
have made at least one more move because the final move to reach
the diagram must have been Nc3-d5.  It turns out that White has
just one way to play exactly one extra move: instead of Bb2,
play Ba3,Bb2.\footnote{
  This kind of ``tempo loss'' is a common feature in composed help-games.
  It can feel paradoxical that a player (here White) rushing to reach
  the diagram as quickly as possible must deliberately waste moves,
  and even more surprising when there is just one viable way
  to waste the right number of moves.
  }
So the minimal games have Black playing the six moves above,
and White playing the \hbox{three-move} sequence b3,Ba3,Bb2,
the \hbox{two-move} sequence Nf3,Nh4, and the single move Nc3
in some order, and then checkmating with Nd5.  Hence the number
of minimal games is ${6 \choose {3,2,1}} = 60$ as desired.
Once again this number is not reduced by accidental checks,
because Black plays d6 first and Kf6 last, so cannot be prematurely
checked by Ba3 or Bb2.

{\em Exercise}\/:  By the time this article appears,
Richard Stanley will be $61$ years old.
Construct an appropriate enumerative chess problem 
as a birthday tribute.\footnote{
  I'm told that in some traditions it is one's 61st birthday,
  not the 60th, that is regarded as an important milestone.
  }
You may, but do not have to, exploit the fact that $61 = E_6$.


\vspace*{2ex}

\begin{thebibliography}{EKLP}
\bibitem[CEF]{CEF} Chow, T.Y., Eriksson, H., Fan, C.K.:
  Chess Tableaux.
  {\em Electronic J.\ Combinatorics} {\bf 11} (2), 2004--2005
  (published June 14, 2005).
  Online at
  {\sf http://www.combinatorics.org/Volume\_11/Abstracts/v11i2a3.html}
\bibitem[El1]{NDE:zeta} Elkies, N.D.:
  On the Sums $\sum_{k=-\infty}^\infty (4k+1)^{-n}$,
  {\em Amer.\ Math.\ Monthly} {\bf 110} \#7
  \hbox{(8--9/2003), 561--573.
  Nearly isomorphic with {\sf http://arxiv.org/math/0101168}$\;$;}
  Corrigenda: {\em Amer.\ Math.\ Monthly} {\bf 111} \#5 (May 2004),~456.
\bibitem[El2]{NDE:2004} Elkies, N.D.:
  Freshman Seminar 23j: Chess and Mathematics,
  Supplementary Question: An enumerative chess problem (2004).
  Online at
  {\sf http://abel.math.harvard.edu/$\sim$elkies/FS23j.04/init.html}~.
\bibitem[EKLP]{KELP}
  Elkies, N., Kuperberg, G., Larsen, M., Propp, J.:
  {\em Alternating sign matrices and domino tilings},
  {\em J.\ Algebraic Combin.}\ {\bf 1} (1992), 111--132 and 219--234.
\bibitem[FH]{FH} Fulton, W., Harris, J.:
  {\em Representation Theory: A First Course.}
  New York: Springer, 1991 (GTM {\bf 129}).
\bibitem[Puu]{Puusa} Puusa, A.:
  {\em Queue Problems.}
  Finnish Chess Problem Society (Suomen Teht\"av\"aniekat), 1992.
\bibitem[Slo]{Sloane:OEIS} Sloane, N.J.A.:
  {\em The On-Line Encyclopedia of Integer Sequences},
  {\sf http://www.research.att.com/$\sim$njas/sequences}~.
\bibitem[St1]{Stanley:EC1} Stanley, R.P.:
  {\em Enumerative Combinatorics}, Vol.~1.
  New York/Cambridge: Cambridge Univ.\ Press, 1997.
\bibitem[St2]{Stanley:EC2} Stanley, R.P.:
  {\em Enumerative Combinatorics}, Vol.~2.
  New York/Cambridge: Cambridge Univ.\ Press, 1999.
\bibitem[St3]{Stanley:CAT} Stanley, R.P.:
  Exercises on Catalan and Related Numbers,
  excerpted from {\em Enumerative Combinatorics}, Vol.~2,
  23 June 1998, online at
  {\sf http://www-math.mit.edu/$\sim$rstan/ec/catalan.pdf};
  Solutions to Exercises on Catalan and Related Numbers,
  online at {\sf http://www-math.mit.edu/$\sim$rstan/ec/catsol.pdf}~.
  See also: Catalan Addendum, 27 May 2005, online at
  {\sf http://www-math.mit.edu/$\sim$rstan/ec/catadd.pdf}~.
\bibitem[St4]{Stanley:ST} Stanley, R.P.:
  Queue problems revisited,
  to appear in {\em Suomen Teht\"av\"aniekat}.
\bibitem[WF]{SPG} Wilts, G., Frolkin, A.:
  {\em Shortest Proof Games}.
  Karlsruhe: Gerd Wilts, 1991.
\end{thebibliography}
\end{document}